\magnification =\magstep1
\baselineskip=14pt
\overfullrule=0pt

\font\footfont=cmr8
\font\footitfont=cmti8

\font\largei=cmti13 at 15pt
\font\medi=cmti12 at 12 pt
\font\bb=msbm10
\font\sbb=msbm8
\font\abs=cmr9

\def\dim{\hbox{dim}\hskip 0.3em}

\def\bbC{\hbox{\bb C}}
\def\sbbC{\hbox{\sbb C}}
\def\qed{\hfill\hbox{\vrule\vtop{\vbox{\hrule\kern1.7ex\hbox{\kern1ex}%
  \hrule}}\vrule}\vskip1.5ex}
\def\eqed{\hbox{\quad
  \hbox{\vrule\vtop{\vbox{\hrule\kern1.7ex\hbox{\kern1ex}\hrule}}\vrule}}}

\def\margin#1#2{\rlap{\hbox to #2em{\ \hfill}\hskip-2em*#1*}}

\def\obeyspaces{\catcode`\ =\active\catcode`\	=\active}
{\obeyspaces\gdef {{\hskip.5em}}\gdef	{{\hskip2em}}}
{\catcode`\^^M=\active %
	\gdef\obeylines{\catcode`\^^M=\active \gdef^^M{\vskip.0ex\ }}}%

\def\idMaccode{\vskip3ex\bgroup\parindent=0pt\baselineskip=11pt
	\obeylines\obeyspaces\bf}
\def\adMaccode{\catcode`\#=11\catcode`\%=11\baselineskip=9pt\tt}
\def\dMaccode{\vskip2ex\bgroup\parindent=0pt%
	\baselineskip=9pt
	\obeylines
	\obeyspaces
	\catcode`\^=11\catcode`\_=11\catcode`\#=11\catcode`\%=11\catcode`\{=11\catcode`\}=11
	\tt}
\def\enddMaccode{\unskip \egroup\vskip-1ex}
\def\Maccode{\bgroup
\catcode`\^=11\catcode`\_=11\catcode`\#=11\catcode`\%=11\catcode`\{=11\catcode`\}=11 \tt"}
\def\endMaccode{\unskip" \egroup}

\def\label#1{{\global\edef#1{\the\sectno.\the\thmno}}}
\def\pagelabel#1{{\global\edef#1{\the\pageno}}}

\def\isnameuse#1{\csname #1\endcsname}

\def\issecond#1#2{#2}
\def\isifundefined#1#2#3{
	\expandafter\ifx\csname #1\endcsname\relax #2
	\else #3 \fi}
\def\pageref#1{\isifundefined{is#1}
	{{\bf ??}\message{Reference `#1' on page [\number\count0] undefined}}
	{\edef\istempa{\isnameuse{is#1}}\expandafter\issecond\istempa\relax}}

\def\today{\ifcase\month\or January\or February\or March\or
  April\or May\or June\or July\or August\or September\or
  October\or November\or December\fi
  \space\number\day, \number\year}

\newcount\sectno \sectno=0
\newcount\thmno \thmno=0
\def \section#1{\vskip 0.3truecm\relax
	\global\advance\sectno by 1 \global\thmno=0
	\noindent{\bf \the\sectno. #1} \vskip 0.2truecm\relax}
\def \subsection#1{\vskip 0.5truecm\relax
	\noindent{\bf #1} \vskip 0.4truecm\relax}
\def \thmline#1{\vskip .2cm\relax
	\global\advance\thmno by 1
	\noindent{\bf #1\ \the\sectno.\the\thmno:}\ \ %
	\bgroup \advance\baselineskip by -1pt \it
	\abovedisplayskip =4pt
	\belowdisplayskip =3pt
	\parskip=0pt
	}
\def \wthmline#1{\vskip .2cm\relax
	\noindent{\bf #1:}\ \ %
	\bgroup \advance\baselineskip by -1pt \it
	\abovedisplayskip =4pt
	\belowdisplayskip =3pt
	\parskip=0pt
	}
\def \dthmline#1{\vskip 6pt\relax
	\noindent{\bf #1:}\ \ }
\def \defin{\thmline{Definition}}
\def \thm{\thmline{Theorem}}

\def \endb{\egroup \vskip 0.2cm\relax}

\def \lemma{\thmline{Lemma}}
\def \remark{\dthmline{Remark}}

\def \example{\dthmline{Example}}

\def \proof{\vskip0.0cm\relax\noindent {\sl Proof:\ \ }}

\ %
\vskip 3ex
\centerline{\largei  Computing Instanton Numbers of Curve Singularities}

\vskip 4ex
\centerline{\medi Elizabeth Gasparim
 and Irena Swanson\unskip\footnote{*}{
\footfont The first author  acknowledges support from NSF NMSU ADVANCE
grant NSF-0123690 and the second author  acknowledges
support from NSF grant DMS-9970566.}}

\unskip\footnote{ }{{\footitfont 2000 Mathematics Subject Classification}
\footfont 13D99, 13P99.}


\bigskip

\bgroup
\narrower\narrower
\baselineskip=11pt
\noindent {\abs
We present an algorithm for computing instanton numbers of curve singularities.
A comparison is made between these
and some other  invariants of curve singularities.
The algorithm is implemented
in  Macaulay2,
and can be downloaded from
 http://www.math.nmsu.edu/\~{}iswanson/instanton.m2.}

\egroup
\bigskip
\medskip

We start with any polynomial $p(x,y)$ defining a plane
curve with singularity
at $ 0 \in {\bbC}^2$.
Let $\pi\colon \widetilde{{\bbC}^2} \rightarrow {\bbC}^2$ denote the blow-up
of ${\bbC}^2$ at the origin and let $j$ be a positive integer.
The data $(j,p)$  determines
a holomorphic bundle $E(j,p)$
on $\widetilde{{\bbC}^2}$ with splitting type $j$
and extension class $p$.
We algorithmically compute numerical invariants of the
bundle $E(j,p)$ and use them as invariants of the curve.

Tables giving examples of 
  these new invariants of the curve
together with some classical invariants
are given  in Section~1. A specially interesting example appears on 
table IV, where we give two inequivalent singularities that have 
all the same classical invariants, but are distinguished by instanton 
numbers.

The holomorphic bundle $E(j,p)$ and its numerical invariants
have interpretation in mathematical physics as instantons
on $\widetilde{{\bbC}^2}$ and  numerical invariants of the instantons.
Here we mention  briefly some properties of these invariants.
More details are given in [G3].
Instantons are well known to have a topological
invariant called the charge, which in the case of compact surfaces
corresponds to a second Chern number of a vector bundle.
For the case of bundles on a blown-up surface,
the local Chern number around the exceptional divisor
decomposes as a sum of two numerical invariants (denoted $w$ and $h)$,
which are local analytic invariants of the bundle.
In [BG] it is shown that the pair of invariants $(w,h)$
is strictly finer than the Chern numbers in the following sense.
The pair $(w,h)$ provides the coarsest
stratification of moduli of instantons on the blown-up plane for which the
strata are Hausdorff.
In contrast,
the stratification by topological charge does not provide Hausdorff strata
[BG, Theorem 4.1].
Hence, the pair $(w,h)$
gives strong numerical invariants, detecting more than topological information.

The idea of using these invariants for curves is natural given that,
as a first step toward resolving the singularity of the plane curve defined 
by $p,$
one blows up the plane,  thus arriving at $\widetilde{{\bbC}^2}$,
the base space of the bundles we construct.

Firstly, let us explain how to pass from curves on 
${\bbC}^2$ to vector bundles on 
$\widetilde{{\bbC}^2}.$
In what follows
we fix, once and for all, the following coordinate charts.
Let $(x,y)$ be the coordinates on ${\bbC}^2.$ Write
$\widetilde{{\bbC}^2} = U \cup V$ where
$$U = \{(z,u)\} \cong {\bbC}^2 \cong \{(\xi,v)\} = V$$
with $(\xi,v)=(z^{-1}, zu)$
in $U \cap V$.
Note that in these coordinates, the blow-up map
$\pi\colon  \widetilde{{\bbC}^2} \rightarrow {{\bbC}^2}$
gives 
$x = u$, $y = zu.$ We denote by $\ell $ the exceptional divisor.

Once these charts are fixed, there is a
canonical choice of transition
matrix for bundles on $\widetilde{{\bbC}^2}.$
If $E$ is a holomorphic rank 2 bundle on $\widetilde{{\bbC}^2}$
with $c_1(E)=0$ then there exists an integer $j$,
called the {\it splitting type} of $E$,
such that  $E|_{\ell} \cong {\cal O}(j) \oplus {\cal O}(-j)$.
By [G2, Theorem~2.1]
$E$ has a canonical transition matrix
of the form
$$ \left(\matrix{ z^j & q \cr 0 & z^{-j} \cr} \right)
\eqno (1)
$$
from $U$ to $V$,
where $$ q \colon = \sum_{i = 1}^{2j-2} \sum_{l = i-j+1}^{j-1}q_{il}z^lu^i$$
is a polynomial in $z, \,z^{-1}$ and $ u.$

In this paper, we start with a polynomial $p \in {\bbC} [x,y]$ 
and for each positive integer $j$ we associate to the pair $(j,p)$ 
a holomorphic bundle $E(j,p)$ on $\widetilde{{\bbC}^2}$ 
obtained as follows.


\defin
Given  a polynomial $p(x,y)$ defined over  ${\bbC}[x,y]$ we  first 
consider the polynomial 
$p(u,zu) = \sum_i\sum_l p_{il}z^lu^i$ 
 obtained from $p$ 
by making $x=u$ and $y=zu.$ 
If furthermore a positive  integer $j$ is chosen, we define 
$\bar p =\sum_i\sum_l \bar p_{il}z^lu^i$  by the rule

$$ \bar p_{il} = \left\{ 
\matrix {  p_{il} & {if}\,\,   1 \leq i \leq 2j-2 \,\,{and} 
\,\, 0 \leq l\leq  j-1; \hfil \cr
  0  &  otherwise. \hfil }\right.\eqno(2)$$
We then set 
$E:=E(j,\bar p)$ to be the bundle
given in canonical coordinates (as above)  by the transition matrix
$$
T = \left(\matrix{ z^j & \bar p \cr 0 & z^{-j} \cr} \right).
\eqno (3)
$$
\endb

We have therefore associated to each pair $(j,p)$ formed by an integer and 
a polynomial $p(x,y)$ a unique vector bundle 
$E(j,\bar p)$ over $\widetilde{{\bbC}^2}.$ We now define 
 instanton invariants
$w$ and $h$ for these vector bundles.

\defin
\label{\defQRone}
Given a bundle $E$ over
$\widetilde{{\bbC}^2}$ we define the sheaf $Q$ by the exact sequence
$$0 \rightarrow \pi_*E \rightarrow (\pi_*E)^{\vee\vee} \rightarrow
Q \rightarrow 0$$ and set
$$w(E)\colon = l(Q), \,\,\,\,\,\,\,\,
h(E)\colon = l(R^1\pi_*E).$$
\endb

In instanton terminology, $w$ is called the width of the instanton
and $h$ is called the height of the instanton.
The charge of the instanton is $c=w+h$.
These are clearly analytic invariants of the bundle.
To use them as invariants of the singularity,
and to perform the calculations,
we choose charts for the bundle and a representative $p$ of the singularity.
Such choice works  particularly well in the case of quasi-homogeneous 
singularities, where preferred representatives are known to exist,
cf. Saito [Sa] or Arnold [Ar, page 95].

If a different representative $p'$
of the germ of the singularity is chosen, then
by classical theory it is known that there exists a
holomorphic change of coordinates taking $p$ to $p'$.
The same change of coordinates has to be made to take the bundle
$E(j,\bar p)$ to a bundle $ E'(j,p')$.
$E$ and $ E'$ are isomorphic bundles
and therefore have the same invariants.
In this sense the numbers can be considered as
analytic invariants for germs of curve singularities.

In this paper we prove that the computation of $h$ and $w$ can be automated.
We implemented our algorithm in the symbolic computer algebra package
Macaulay2 (due to Grayson-Stillman [GS]).
Having implemented  algorithms to compute both invariants,
we obtained a large amount of examples,
which lead us to 
a simple formula  for $h.$ 
\medskip

\noindent {\bf Theorem 3.3:} 
{\it Let $m$ denote the largest power of $u$ dividing $\bar p$.
If $E$ is the bundle defined by data $(j,\bar p)$,
then
$
h(E) = {j \choose 2} - {j-m \choose 2}.
$}
\medskip

The algorithm we present here has as 
 input data the polynomial $p(x,y)$ and the integer 
$j$, and as output the values of 
the instanton height $h(E(j,\bar p))$  and width
$w(E(j,\bar p))$. 
Our implementation on Macaulay2 is for $p$ with rational coefficients.
%

\section{Some tables}

Before giving the algorithm, we illustrate some 
instanton numbers by tabulating them together with
some  classical invariants.
Tables I and II give examples of instanton numbers of monomials.

\def\vertline{{\vrule height2.7ex depth.8ex}}
\def\chalign#1{{$$\vbox{\halign{#1}}$$}}

\chalign{
\vertline\ \hfil{\it #}\hfil\vadjust{\hrule\vskip-2pt} \vertline &&
	\hskip0.5em\hfil{\it #}\hfil\hskip0.5em \vertline \cr
\noalign{\hrule}
TABLE I & \multispan3 \hfil \it splitting type $2$ \hfil \vertline \cr
\noalign{ \vskip 1pt}
polynomial & $w$ & $ h$ & charge \cr
$x$ & 1 & 1 & 2 \cr
$y$ & 1 & 1 & 2 \cr
$xy$ & 2 & 1 & 3 \cr
}

\chalign{
\vertline\ \hfil{\it #}\hfil\vadjust{\hrule\vskip-2pt} \vertline &&
	\hskip0.5em\hfil{\it #}\hfil\hskip0.5em \vertline \cr
\noalign{\hrule}
TABLE II & \multispan3 \hfil \it splitting type $3$ \hfil \vertline \cr
\noalign{ \vskip 1pt}
polynomial & $w$ & $h$ & charge \cr
$x$ & 1 & 2 & 3 \cr
$y$ & 1 & 2 & 3 \cr
$x^2$ & 3 & 3 & 6\cr
$xy$ & 2 & 3 & 5 \cr
$y^2$ & 3 & 3 & 6 \cr
$x^2y$ & 4 & 3 & 7\cr
$xy^2$ & 4 & 3 & 7\cr
$x^2y^2$ & 5 & 3 & 8\cr
}

\noindent Table III 
illustrates that for a  fixed polynomial,
we may get different invariants as we vary the splitting type.

\chalign{
\vertline\ \hfil{\it #}\hfil\vadjust{\hrule\vskip-2pt} \vertline &&
\hskip0.5em\hfil{\it #}\hfil\hskip0.5em \vertline \cr
\noalign{\hrule}
TABLE III & \multispan3 \hfil \it splitting type $7$ \hfil \vertline
& \multispan3 \hfil \it splitting type $8$ \hfil \vertline \cr
\noalign{ \vskip 1pt}
polynomial & $w$ & $h$ & charge & $w$ & $h$ & charge \cr
$x^2-y^2$ & 2 & 11 & 13 & 3 & 13 & 16 \cr
$x^2-y^3$ & 3 & 11 & 14 & 3 & 13 & 16\cr
$x^2-y^5$ & 3 & 11 & 14 & 3 & 13 & 16\cr
$x^3-y^3$ & 4 & 15 & 19 & 4 & 18 & 22\cr
$x^3-y^4$ & 6 & 15 & 21 & 6 & 18 & 24\cr
$x^3-y^5$ & 6 & 15 & 21 & 6 & 18 & 24\cr
}

Our interest in instanton numbers was partially fueled by the fact
that in some cases,
instanton numbers give finer information than the classical invariants.
We considered the invariants:
multiplicity $m$,
$\delta_P= \dim(\widetilde{\cal O}/{\cal O})$,
$\mu=$ Milnor number,
and $\tau=$ Tjurina number.
Note that the Milnor and Tjurina numbers 
are defined only for isolated singularities,
but instanton numbers
are well defined for non-isolated singularities as well.
Here is an example where  instanton numbers distinguish the
singularities,
but other invariants do not.

\chalign{
\vertline\ \hfil{\it #}\hfil\vadjust{\hrule\vskip-2pt} \vertline &&
\hskip0.5em\hfil{\it #}\hfil\hskip0.5em \vertline \cr
\noalign{\hrule}
\multispan5\vertline\ TABLE IV \hfil \vertline&\multispan3 \it splitting type $3$ \hfil \vertline \cr
\noalign{\vskip -1pt}
\noalign{\hrule}
polynomial & m & $\delta_P$ & $\mu$ & $\tau$ & $w$ & $h$ & charge \cr
$x^3-x^2y+y^3$  & 3&  3 & 4 & 4 & 4 & 3 & 7\cr
$x^3-x^2y^2+y^3$ & 3& 3 & 4 & 4 & 5 & 3 & 8\cr
}

\noindent
Note that the first polynomial is irreducible, whereas
the second is reducible in the local ring, cf. [Ha, ex. I 5.14],
so they define inequivalent singularities.

We believe that the 
fact that the classical invariants $m$, $\delta_P, \mu$, and $\tau$ 
do  not distinguish some inequivalent singularities 
is an evidence that finer invariants are useful.
In tables V-VII below we give examples of inequivalent singularities
which are distinguished by instanton numbers,
and also by one other classical invariant.

\chalign{
\vertline\ \hfil{\it #}\hfil\vadjust{\hrule\vskip-2pt} \vertline &&
\hskip0.5em\hfil{\it #}\hfil\hskip0.5em \vertline \cr
\noalign{\hrule}
\multispan5\vertline\ TABLE V \hfil \vertline&\multispan3 \it splitting type $4$ \hfil \vertline \cr
\noalign{\vskip -1pt}
\noalign{\hrule}
polynomial & $m$ & $\delta_P$ & $\mu$ & $\tau$ & $w$ & $h$ & charge \cr
$x^2-y^7$ & 2 & 3 & 6 & 6 & 3 & 5 & 8\cr
$x^3-y^4$ & 3 & 3 & 6 & 6 & 6 & 6 & 12\cr
}

\noindent Table VI
comes from [Ha, ex. V 3.8] on page 395.
However,
in the statement of this exercise,
the first polynomial contains an incorrect exponent:
it is written as ``$x^4y - y^4$" but it should be ``$x^5y - y^4$".

\chalign{
\vertline\ \hfil{\it #}\hfil\vadjust{\hrule\vskip-2pt} \vertline &&
\hskip0.5em\hfil{\it #}\hfil\hskip0.5em \vertline \cr
\noalign{\hrule}
\multispan5\vertline\ TABLE VI \hfil \vertline&\multispan3 \it splitting type $4$ \hfil \vertline \cr
\noalign{\vskip -1pt}
\noalign{\hrule}
polynomial & $m$ & $\delta_P$ & $\mu$ & $\tau$ & $w$ & $h$ & charge \cr
$x^4-xy^5$ & 4 & 9 & 17 & 17 & 10 & 6 & 16\cr
$x^4-x^2y^3-x^2y^5-y^8$ & 4 & 9 & 17 & 15 & 8 & 6 & 14\cr
}

%

\noindent Table VII comes from  the list of bimodular singularities
given by Arnold in~[Ar, page 159].

\chalign{
\vertline\ \hfil{\it #}\hfil\vadjust{\hrule\vskip-2pt} \vertline &&
\hskip0.5em\hfil{\it #}\hfil\hskip0.5em \vertline \cr
\noalign{\hrule}
\multispan5\vertline\ TABLE VII \hfil \vertline&\multispan3 \it
 splitting type $8$ \hfil \vertline \cr
\noalign{\vskip -1pt}
\noalign{\hrule}
polynomial & $m$ &  $\delta_P$ & $\mu$ & $\tau$ & $w$ & $h$ & charge \cr
$x^3+x^2y^3+y^9+xy^7$ & 3 & 9 & 15 & 16 & 6 & 18 & 24\cr
$x^3y+x^2y^3+xy^6+y^7$ & 4 & 9 & 14 & 15 & 7 & 22 & 29\cr
$x^4+x^2y^3+y^6$ & 4 & 8 & 15 & 15 & 10 & 22 & 32\cr
$(x^2+y^3)^2+xy^4$ & 4 & 7 & 12 & 13 & 9 & 22 & 31\cr
$(x^2+y^3)^2+xy^3$ & 4 & 6 & 9 & 9 & 6 & 22 & 28\cr
}

\noindent Table VIII below comes from the list of 
exceptional families of unimodal singularities in ~[Ar, page 95].


\chalign{
\vertline\ \hfil{\it #}\hfil\vadjust{\hrule\vskip-2pt} \vertline &&
\hskip0.5em\hfil{\it #}\hfil\hskip0.5em \vertline \cr
\noalign{\hrule}
TABLE VIII \hfill & \multispan3 \it splitting type $7$ \hfil \vertline \cr
\noalign{ \vskip 1pt}
polynomial & $w$ & $h$ & charge \cr
$ x^3+y^7+xy^5$ & 6 &15 & 21\cr
$ x^3+y^8+xy^6$ & 6 &15 & 21\cr
$x^3y+xy^4+x^2y^3$ & 7 & 18 & 25\cr
$x^3+xy^5+y^8$ &  6& 15  & 21\cr
$x^3y+y^5+xy^4$& 7& 18& 25\cr
$x^3y+y^6+xy^5$& 7& 18  & 25 \cr
$x^4+xy^4+y^6$ & 9 & 18 & 25\cr
}

\section{Computing the instanton width}

In this section we present an algorithm that takes as input the  
pair $(j,p)$  and computes
the instanton width 
 $w=l(Q)$,
where $Q$ is the skyscraper sheaf as given in  Definition~\defQRone.
Since $Q$ is supported at zero,
$l(Q)$ is the dimension of the stalk at zero.
Let $M$ be the completion of the stalk $\pi_*E$ at zero,
that is,
$M \colon = (\pi_*E)_0^{\wedge}$.
Then the length of $Q$ 
equals the dimension of the cokernel of the canonical map from $M$
to its double dual.
If we can compute $M$,
we can also compute $Q$ via the following lemma:

\lemma
\label{\lmQ}
Let $R$ be a commutative Noetherian ring
and $A$ an $n \times m$ matrix with entries in $R$.
Let $M$ be the $R$-module such that
$R^m {\buildrel A \over \longrightarrow} R^n \to M \to 0$ is an exact sequence.
Let $N$ be the kernel of the transpose of $A$.
Then $N$ is a submodule of $R^n$,
say generated by $t$ elements.
Let $B$ be the $n \times t$ matrix
whose columns are the given generators of $N$.
Let $C$ be the matrix
such that $R^k {\buildrel C \over \longrightarrow} R^t \to N \to 0$ is exact.
Let $\underline{\hbox{\ }}^{\vee} = Hom_R(\underline{\hbox{\ }}\ , R)$.
Then
\item{(i)}
$M^{\vee\vee}$ is isomorphic to the kernel of the transpose of $C$,
\item{(ii)}
the image of the canonical map $M \to M^{\vee\vee}$
is isomorphic to the module generated by the columns of the transpose of $B$,
\item{(iii)}
$M^{\vee\vee}/ M$ is isomorphic to $(kernel\ C^T)/(image\ B^T)$.
\endb

\proof
By left-exactness of the $Hom$ functor,
$M^{\vee} = Hom_R(M,R)$
is the kernel of the map defined by the transpose of $A$.
Thus $N = M^{\vee}$.
Similarly,
by the definition of $C$,
$M^{\vee\vee}$ is the kernel of the transpose of $C$.
Let $f_1, \ldots, f_t$ generate $M^{\vee} \subseteq R^n$.
The natural map $M \to M^{\vee\vee}$ takes $m$ to $\varphi_m$,
where for each $f \in M^{\vee}$,
$\varphi_m(f) = f(m)$.
Thus, in coordinates,
the image of $m \in M$ in $M^{\vee\vee} \subseteq R^t$
is the vector $(f_1(m), \ldots, f_t(m))$,
therefore the image of $M$ in $M^{\vee\vee}$
is generated by the columns of the transpose of $B$.
\qed

Once $Q$ is computed as in the previous lemma,
its length can be computed as well.
Thus it remains to compute the $\bbC[[x,y]]$-module $M$.
By the Theorem on Formal Functions,
$$M \cong \lim_{\longleftarrow}
H^0(\ell_n, E|{\ell_n})$$
as ${\bbC}[[x,y]]$-modules,
where $\ell_n $ is the $n$th formal neighborhood of $\ell$.
We use the following lemma:
\lemma [G1, Lemma 2.2]
\label{\width}
 Set $\bar M= H^0(\ell_{2j-2},{E}|_{\ell_{2j-2}}),$
 and let $\rho \colon \bar M \hookrightarrow 
\bar M^{\vee\vee}$ be the natural inclusion of $\bar M $ into its bidual.
Then $l(Q)=$ dim coker$\,\rho.$
\endb

This lemma greatly simplifies the calculations. Using this lemma,
for the purpose of finding $l(Q),$ we may assume $M=\bar M.$ 
Under this identification,
and  with our choice of transition matrix $T$ in $(3),$ it then follows that 
$M$ is generated by all holomorphic vectors $v$ over $\bbC[[z,u]]$
for which $T v$ is a holomorphic vector over $\bbC[[z^{-1},zu]]$.
Thus, we need to find all vectors $v= (a, b)$,
given by $a = \sum_{i,l} a_{il} z^l u^i $,
$b = \sum_{i,l} b_{il} z^l u^i $,
where all $i$ and $l$ are non-negative,
and $a_{il}, b_{il}$ are (unknowns) in the field,
for which $Tv$ is holomorphic over $\bbC[[z^{-1},zu]]$.
This restriction yields relations among the unknowns $a_{il}, b_{il}$.
To get all such relations,
at the beginning we treat $a_{il}$ and $ b_{il}$ as variables.


\lemma
Whenever $l > i + j$,
then $b_{il} = 0$.
Whenever $l > i$,
then $a_{il} = 0$.
\endb

\proof
The second coordinate of $T v$ is
$z^{-j} b = \sum_{i,l} b_{il} z^{l-j}u^i
= \sum_{i,l} b_{il} z^{l-j-i}(zu)^i$.
In order for this to be holomorphic in $z^{-1}$ and $zu$,
necessarily the coefficients $b_{il} $ with $l > i + j$ must vanish.

The first entry of $T v$ is
$$
\eqalignno{
z^j \sum_{i,l} a_{il} z^l u^i &+
\bar p \sum_{i_1 ,l_1 } b_{i_1 l_1 } u^{i_1} z^{l_1} \cr
&
= \sum_{i,l} a_{il} z^{l-i+j}(zu)^i
+ \sum_{i_0 = 1}^{2j-2} \sum_{l_0 = i_0-j+1}^{j-1} \bar p_{i_0l_0}(zu)^{i_0}z^{l_0-i_0}
\sum_{i_1,l_1} b_{i_1l_1} z^{l_1-i_1}(zu)^{i_1}.
}
$$
Let $m$ be the minimum $u$-degree occurring in $\bar p$,
i.e., the minimum $i_0$ such that $\bar p_{i_0l_0} \not = 0$.
Then for each $i < m$,
the coefficient of $u^i$ in the first entry of $T v$
is $\sum_{l} a_{il} (zu)^i z^{l-i+j}$.
In order for this to be holomorphic in $zu$ and $z^{-1}$,
necessarily for all $l > i-j$, $a_{il} = 0$.
Now consider the case $i \ge m$.
The coefficient of $(zu)^i z^{l+j-i}$ in the first entry of $T v$ is
$$
a_{il} +
\sum_{i_0 = 1}^{2j-2} \sum_{l_0 = i_0-j+1}^{j-1} \bar p_{i_0l_0} b_{i-i_0,l+j-l_0},
$$
where $b_{i',l'}$ is treated as zero whenever $i'$ or $l'$ is negative.
By the established bounds,
whenever $l > i-i_0+l_0$,
then $b_{i-i_0,l+j-l_0} = 0$.
Since $\bar p$ is a subpolynomial of $p(u,zu)$,
the only pairs $(i_0,l_0)$ to consider are those with $l_0 \le i_0$.
Thus whenever $l > i-i_0+l_0$,
then $b_{i-i_0,l+j-l_0} = 0$.
In particular,
if $l > i$,
then all $b_{i-i_0,l+j-l_0}$ are zero,
so that the coefficient of $(zu)^i z^{l+j-i}$ in the first entry of $T v$ is
$a_{il}$.
But if $Tv$ is to be holomorphic in $zu, z^{-1}$,
then necessarily $a_{il} = 0$.
\qed

Thus for each $i$,
we need only to consider finitely many unknowns $a_{il}, b_{il}$
to construct the vectors $v$ as in the set-up above,
and find relations on these.
Using Lemma~\width,  we may assume that 
$M= H^0(\ell_{2j-2},{E}|_{\ell_{2j-2}})$.
This means that
we only need relations involving $a_{il}, b_{il}$ for $i \le 2j-2$.
But the relations involving $b_{il}$ with $i \le 2j-2$
arising from $z^j a + \bar p b$ being holomorphic over $\bbC[[u,zu]]$
involve variables $a_{i'l}$ and $b_{i'l}$
with $i' \le 2j-2 + deg_u \bar p$,
so that

\lemma
\label{\lmjbd}
Let $N$ be the sum of $2j-2$ plus the $u$-degree of $\bar p$.
Then $M$ is generated by vectors $(a,b)$
with $a = \sum_{i,l} a_{il} z^l u^i $,
$b = \sum_{i,l} b_{il} z^l u^i $,
where all $i$ and $l$ are non-negative,
$a_{il}, b_{il} \in \bbC$,
and
\item{(i)}
$a_{il} = b_{il} = 0$ whenever $i > N$.
\item{(ii)}
$b_{il} = 0$ whenever $l > i + j$.
\item{(iii)}
$a_{il} = 0$ whenever $l > i$.
\qed
\endb

The conditions on the $a_{il}$ and $b_{il}$ as in the lemma above
are not the only ones
that arise from the condition that $Tv$ be holomorphic in $z^{-1}, zu$.
Finding all the conditions amounts to finding the generators of $M$.

\defin
All the remaining relations on the variables $a_{il}, b_{il}$
arise from the condition that in the first entry of $Tv$,
whenever $l > i$,
then the coefficient of $u^i z^l$ must be $0$.
From now on,
we refer to these coefficients as the {\bf generating relations}.
\endb

Note that these coefficients
are all linear forms in the ring $\bbC[a_{il}, b_{il}]$.

As in [G1],
we find these relations successively in the zeroth through $N$th neighborhoods.
Below follows our algorithm which computes the generating relations,
arising from ensuring that the first coordinate of $Tv$ is holomorphic
in $zu, z^{-1}$.
We write all algorithms in this paper in pseudocode,
close to the Macaulay2 code that we implemented.

\idMaccode
Algorithm ``getrelations" to get the generating relations on the $\bf a_{il}, b_{il}$
Input: non-negative integer N,
\hphantom{Input:} fTv = first entry of Tv
\hphantom{Input:} ring R, the polynomial ring in u, z, all$\bf a_{il}, b_{il}$
Output: ideals nonfree and relations in R
\adMaccode
   nonfree = zero ideal in R
   relations = zero ideal in R
   k = 0
   while (k $\le$ N) (
      tempoly = truncation of fTv to terms of $u$-degree at most k
      while (tempoly != 0) (
         tempterm = leading term of tempoly
         i = u exponent of tempterm
         l = z exponent of tempterm
         partp = coefficient of u$^i$z$^l$ in tempoly, linear form in a$_{il}$, b$_{il}$
         tempoly = tempoly - partp * u$^i$z$^l$
         fTv = fTv - partp * u$^i$z$^l$
         if (l > i) then (
            relations = relations + ideal(partp)
            nonfree = nonfree + ideal (leading variable in partp)
         )
      )
      k = k + 1
   )
   return {nonfree, relations}
)
\enddMaccode

The ideal of all the generating relations
obtained in this way is called {\tt relations}.
The ideal {\tt nonfree} contains the leading variable of each relation:
in the sense of linear algebra,
these leading variables are not free.
We use the following ordering of the monomials.

\defin
\label{\deford}
We say that $a_{il} > a_{i'l'}$
if $i > i'$ or if $i = i'$ and $l > l'$,
and similarly $b_{il} > b_{i'l'}$
if $i > i'$ or if $i = i'$ and $l > l'$,
and furthermore that $a_{il} > b_{i'l'}$ for all $i, l, i', l'$.
\endb

Thus by the form of $Tv$,
each relation contains at most one $a_{il}$,
and each $a_{il}$ appears at most once in a generating relation.
Observe that if all the coefficients of 
$p(x,y)$ are in a subfield $F$ of $\bbC$,
then all the generating relations have coefficients in $F$.
The algorithm \Maccode getrelations \endMaccode above computes
the relations among the given $a_{il}, b_{il}$,
but some of these relations are ``fake" in the following sense:

\example
Let $p(x,y) = x^2-y^3, j = 3$.
Then with notations as above,
$N = 7$,
$\bar p = u^2 - u^3 z^3$.
For holomorphic $a = \sum a_{il} u^i z^l$, $b = \sum b_{il} u^i z^l$,
the coefficient of $u^8z^9$ in $z^j a + \bar p b$
is $0 = a_{86} + b_{69} - b_{56}$.
However,
if we restrict the first index of $a_{il}$ to only vary from $0$ to $N$,
then \Maccode getrelations \endMaccode gives the ``fake'' relation
$0 = b_{69} - b_{56}$.

Thus, computation of the instanton width
will have to account for and remove such ``fake" relations.
We do this as follows.
Such relations only involve variables $a_{il}, b_{il}$
with $i > 2j-2$.
Thus these variables are not allowed to be free variables
in the sense of linear algebra.
Using Lemma~\width, 
the remaining free variables do give a generating set of $M$ as follows.
For each of the free variables,
set that variable to $1$ and all the others to $0$ in $(a,b)$.
This produces a finite generating set of $M$ as a module over $\bbC[[x,y]]
=\bbC[[u,zu]]$
of elements whose entries are in $\bbC[[u,z]]$.
As $u$ is a non-zerodivisor,
$M$ is isomorphic to $u^k M$ for arbitrary integer $k$.
By Lemma~\lmjbd,
$u^j a$ and $u^j b$ are both polynomials in $u (= x)$ and $zu (= y)$,
so that the generators of $u^j M$
can be written as pairs of polynomials in $\bbC[[x,y]]$.

\idMaccode
Algorithm ``polyconv" to convert $\bbC[u,z]$-polynomials to $\bbC[x,y]$-polynomials
Input: a polynomial f in u and z
Output: a ``truncated" polynomial g(x, y) such that g(u,zu) = f',
\hphantom{Output:} where f' is that part of f for which this can be done
\adMaccode
	g = 0 zero element of \bbC[u,z]
	while (f != 0) (
		lf = leading term of f
		i := the u exponent of lf
		l := the z exponent of lf
		if (l $\le$ i) then
			g = g + y$^l$ * x$^{i-l}$ * leading coefficient of lf
		f = f - lf
		)
	return g
	)
\enddMaccode

\vskip-4ex
\idMaccode
Algorithm ``setvectors" to express generators of $u^{N+j} M$
\hphantom{Algorithm:} as vectors with entries in $\bbC[x,y]$
Input: polynomials $\bf u^j a, u^j b$, named Apoly, Bpoly, respectively,
\hphantom{Input:} lists changeables, allvars
Output: A, a presentation matrix for the $\bbC[[x,y]]$-module $M$
\adMaccode
    Mxy =  zero submodule of $\bbC[x,y]$$^2$
    total = #changeables
    k = 0
    while (k < total) (
        tapoly = substitute in Apoly the k$^{th}$ changeable variable to 1
        tbpoly = substitute in Bpoly the k$^{th}$ changeable variable to 1
        tapoly = substitute in tapoly all other variables to 0
        tbpoly = substitute in tbpoly all other variables to 0
        A = convert tapoly into a polynomial in x and y (use polyconv)
        B = convert tbpoly into a polynomial in x and y (use polyconv)
        Mxy = Mxy + submodule of $\bbC$[x,y]$^2$ generated by (A,B)
        k = k + 1
    )
    return a presentation of the $\bbC[x,y]$-module Mxy
)
\enddMaccode

The output of the last routine
is the presentation matrix of a $\bbC[x,y]$-module,
which by faithful flatness of $\bbC[[x,y]]$ over $\bbC[x,y]$
is also the presentation matrix of the $\bbC[[x,y]]$-module $M$.
Finally,
tying it all together:

\idMaccode
Algorithm to compute instanton width of the instanton with data $\bf (j,p)$
Input: a polynomial p in $\bbC[x,y]$ and a non-negative integer j
Output: the width of the instanton with data $\bf (j,p)$
\adMaccode
    \=p = p(u,zu) truncated to u-degree at most 2j-2
    N = 2j-2 + u-degree of \=p
    R = \bbC[u,z,a$_{il}$, b$_{il}$], i < N+1, ordered as in Definition~\deford
    a = $\bf \sum$ a$_{il}$ u$^i$ z$^l$
    b = $\bf \sum$ b$_{il}$ u$^i$ z$^l$
    fTv = z$^j$ a + \=p b
    compute relations and nonfree variables as in algorithm getrelations
    changeables = all a$_{il}$, b$_{il}$ with $i \le 2j-2$, which are not in nonfree
    a = a after applying all the relations
    b = b after applying all the relations
    A the presenting matrix of the $\bbC$[x,y]-module M,
        output of setvectors(u$^{j+N}$apoly, u$^{j+N}$bpoly, changeables, allvars)
    Q = cokernel of the natural map M $\to$ M$^{\vee\vee}$, as in Lemma \lmQ
    return (length of Q)
)
\enddMaccode

\section{Computing the instanton height}

In this section we compute the instanton height.
Recall that the instanton height $h$ is the length $l(R^1\pi_*E(j,\bar p))$.
Another use of  the Theorem on Formal Functions gives
$$(R^1\pi_* {E})^{\wedge}_0 = \lim_{\longleftarrow}\,\,
H^1(\ell_n, {E}|_{\ell_n}).$$
We use the following lemma:

\lemma  [G1, Lemma 2.3]
\label{\height}
 $ l(R^1\pi_* E) =  dim_{\bbC}\,\,
 H^1(\ell_{2j-2}, E|_{\ell_{2j-2}}).$
\endb
This lemma greatly simplifies the calculations.
We then proceed to compute the first \v{C}ech cohomology.

\remark
We have two options to compute \v{C}ech cohomology.
Given that our
open cover of the base space is given by affine sets that are
open  in the analytic topology as well as in the Zariski topology,
we have the option to compute either 
 holomorphic \v{C}ech cohomology (by taking
holomorphic cochains) or else  algebraic \v{C}ech cohomology (by taking
algebraic, i.e, polynomial cochains). Since $\ell$ is compact,
so are its formal neighborhoods.
By Serre's G.A.G.A.,
holomorphic bundles on a compact variety are algebraic,
therefore
$H^1_{alg}(\ell_n,E\vert_{\ell_n})= H^1_{hol}(\ell_n,E\vert_{\ell_n})$.
It follows that the holomorphic and the algebraic methods give the same answer.
\vskip0.5 ex

We compute the instanton height using holomorphic \v{C}ech cohomology.
The 1-cocycles consist of the vectors $(a,b)$
which are holomorphic functions defined on the intersection $U\cap V$.
Hence $a, b \in \bbC[[u,z,z^{-1}]]$.
The coboundaries consist of the vectors of the form $v +  T^{-1} v'$,
where $v$ is holomorphic in $z,u$ (on $U$) and $v'$ in $z^{-1}, zu$ (on $V$).
First of all we choose simple representatives for the cocycles:

\lemma
\label{\lmcocycleform}
Every 1-cocycle has a representative of the form
$$
\sum_{i=0}^{j-2} \sum_{l=i-j+1}^{-1} \left(\matrix{a_{il}\cr 0}\right) z^lu^i,
$$
with $a_{il} \in \bbC$.
In particular, every 1-cochain represented by 
$ \left(\matrix{a_{il} \cr 0}\right)z^lu^i $ with $i, l \geq 0$ is a coboundary.
\endb

\proof
Let $\sigma$ be a 1-cocycle and let $\sim$ denote cohomological equivalence.
A power series representative for a 1-cochain has the form
$\displaystyle \sigma = \sum_{i = 0}^\infty \sum_{l=-\infty}^{\infty}
\left(\matrix{a_{il} \cr b_{il}}\right) z^l u^i$,
with $a_{il}, b_{il} \in \bbC$.
The 1-cochain
$\displaystyle s_1=
\sum_{i = 0}^\infty \sum_{l=0}^\infty \left(\matrix{a_{il} \cr
b_{il}}\right) z^l u^i$
is holomorphic in $U$,
hence is a coboundary.
Hence $$\sigma \sim \sigma - s_1 =
\sum_{i = 0}^{\infty} \sum_{l=-\infty}^{-1}\left(\matrix{ a_{il} \cr
b_{il}}\right) z^l u^i.
$$
After a change of coordinates
$$T \sigma =\sum_{i = 0}^{\infty} \sum_{l=-\infty}^{-1}
\left( \matrix{z^j a_{il} +
\bar p \, b_{il} \cr
z^{-j} b_{il}}\right) z^l u^i,
$$
but given that
$\displaystyle s_2 =\sum_{i = 0}^{\infty} \sum_{l=-\infty}^{-1}
\left( \matrix{ 0 \cr
z^{-j} b_{il}}\right) z^l u^i$ is holomorphic in $V$,
$$
T \sigma \sim T \sigma -s_2 =\sum_{i = 0}^{\infty} \sum_{l=-\infty}^{-1}
\left( \matrix{z^j a_{il} + \bar p \, b_{il} \cr 0 }\right) z^lu^i$$
and going back to the $U$--coordinate chart,
$$\sigma = T^{-1} T \sigma
\sim \sum_{i = 0}^{\infty} \sum_{l=-\infty}^{-1}
\left( \matrix{ a_{il} +
z^{-j}\bar p \, b_{il} \cr
0 }\right)z^lu^i.
$$
But $\bar p$ contains only terms $z^k$ for $k\leq j-1$,
therefore $z^{-j}\bar p$ contains only negative powers of $z$.
Renaming the coefficients we may write
$\displaystyle \sigma
=\sum_{i = 0}^{\infty} \sum_{l=-\infty}^{-1}
\left( \matrix{ a'_{il}
\cr
0 }\right)z^lu^i$
for some $a'_{il} \in \bbC$,
and consequently $\displaystyle T \sigma
=\sum_{i = 0}^{\infty} \sum_{l=-\infty}^{-1}
\left(\matrix{z^j a'_{il} \cr 0 }\right)z^lu^i$.
Here each term $a'_{il}z^{j+l}u^i$ satisfying
$j+l \leq i$ is holomorphic in the $V$--chart.
Subtracting these holomorphic terms
we are left with an expression for $a$ where the index
$l$ varies as $i-j+1 \le l \le -1$.
This in turn forces $i \leq j-2, $ giving the claimed expression
for the 1-cocycle.
\qed


Thus we only have to consider cocycles of the form $(a,0).$
Which of the cocycles $(a,0)$ is a coboundary?
In other words,
which $(a,0)$ equal to $(c,d) + T^{-1}(c',d')$,
where $c,d$ are holomorphic on $U$ and $c', d'$ on $V$,
or even more simply,
for what $c$ and $d$ holomorphic on $U$ is $T(a+c,d)$ holomorphic on $V$?
The second coordinate of $T(a+c,d)$ is $z^{-j} d$,
and in order for that to be holomorphic,
$d_{il} = 0$ whenever $l \ge i + j$.
This is the only restriction on 
$c$ and $d$ obtained from the second coordinate.
 From  the first coordinate of $T(a+c,d)$ 
we obtain the constraint that $z^j (a+c) + \bar p d$
be holomorphic on $V,$
that is, holomorphic on coordinates $z^{-1} $ and $zu$.

\lemma
\label{\lmnonzerocoh}
Let $E$ be the bundle defined by data $(j,\bar p)$,
and $m$ the smallest exponent of $u$ appearing in $\bar p$.
Then
$$
l(R^1\pi_*E) \geq {j \choose 2} - {j-m \choose 2}.$$
\endb

\proof
Let $\sigma= (a,0)$ denote a 1-cocycle
where  $a= z^lu^i$ 
with $0 \leq i \leq m-1$ and $i-j+1 \leq l \leq -1.$ 
We claim that $\sigma$ 
represents a nonzero cohomology class.  
In fact, for $\sigma$ to be a coboundary there must exist $c$ and $d$,
holomorphic in $U$,
making the expression $z^j (a+c)+\bar p d$ holomorphic in $V$.
However, $z^ja$ is not holomorphic in $V$.
Moreover,
by the choice of $m$,
no term in $\bar p d$ cancels $z^j a$.
Consequently,
no choice of $c$ and $d$ solves the problem of holomorphicity on $V$.
Hence 
$l(R^1\pi_*E)$
is at least the number of independent  cocycles of the form
$\sigma= (a,0)$, 
where $a= z^lu^i$ with  $0 \leq i \leq m-1$ and $i-j+1 \leq l \leq -1.$ 
There are 
${j \choose 2} - {j-m \choose 2}$ such terms.
\qed

\thm
\label{\lmhtatleast}
Let $E$ be the bundle defined by data $(j,\bar p)$,
and $m$ the smallest exponent of $u$ appearing in $\bar p$.
Then
$$
l(R^1\pi_*E) = {j \choose 2} - {j-m \choose 2}.
$$
\endb

\proof
First note that, by Lemma~{\lmcocycleform}, 
if $l \ge j$,
then 
$( z^{l-j} u^i,0)$ is  a coboundary.
Using the proof of Lemma~{\lmnonzerocoh},
it suffices to show that if  $j>l>i \ge m$, then 
$\sigma=( z^{l-j} u^i,0)$ is a coboundary.
Write $\bar p = \sum_{r=0}^s a_r z^{m-r}u^m + p'$,
where $s \in \{0, \ldots, m\}$ is some integer,
$a_r$ are constants,
$a_s \not = 0$,
and $p'$ is a polynomial in $u$ and $zu$
each of whose terms has $u-$degree at least $m + 1$.
Observe that $d=a_s^{-1}u^{i-m} z^{l-m+s}$ is holomorphic on $U,$ since
by assumption $i \ge m$ and $l > m-s.$ 
Therefore $\sigma \sim \sigma'=(z^{l-j}u^i, -d),$
 where $\sim$ denotes 
cohomological equivalence. Changing coordinates we have
$T\sigma'= (z^lu^i- \bar p d, -z^{-j}d)$.
 From $j > l > i \ge m \ge s$ we deduce that $z^{-j}d$ is holomorphic on $V$.
We rewrite the first entry $T\sigma'$ as
$$
 z^l u^i- \bar p d
= \sum_{r = 0}^{s-1} a_s^{-1} a_r  z^{l+s-r}u^i +
a_s^{-1} u^{i-m} z^{l-m+s}p'.
$$
By assumption,  $l-j+s-r \leq s-r \leq s \leq m \leq i,$ therefore
the first sum on the right side of this expression is holomorphic on $V$.
It follows that
$T\sigma'\sim (a_s^{-1} u^{i-m} z^{l-m+s}p', 0).$ 

Let $z^v u^r$ be an arbitrary term in $p'$.
Then $r > m$ and $v \le r$.
If $i-m+r \ge l-m+s+v$,
then the term
$c= a_s^{-1} u^{i-m} z^{l-m+s} \cdot u^r z^v$ is holomorphic on $V$,
therefore $T \sigma'\sim (a_s^{-1} u^{i-m} z^{l-m+s}p'-c , 0).$ 

Removing all such terms $c,$ we may now write 
$T\sigma'\sim (a_s^{-1} u^{i-m} z^{l-m+s}\bar p,0 ),$ where 
$\bar p$  
  contains 
only terms in  $z^v u^r$ such that  $i-m+r < l-m+s+v.$
Consequently $\sigma' \sim  (z^{-j} a_s^{-1} u^{i-m} z^{l-m+s}\bar p,0 )$
and as $i-m+r > i$,
by (reverse) induction on $i$ and $l$, each term 
$( u^{i-m} z^{l-m+s} \cdot u^r z^v z^{-j},0)$ is a coboundary.
Hence, $\sigma'$ is a sum of coboundaries, and is itself a  coboundary.
\qed

Thus when starting with $p(x,y) \in \bbC[x,y]$,
the computation of the instanton height is very fast:
once $m$ is determined,
then the following routine \Maccode iheight \endMaccode
finishes the computation:

\idMaccode
Algorithm ``iheight" to compute instanton height
Input: a polynomial p and a non-negative integer j
Output: returns the length of $\bf R^1$.
\adMaccode
	\=p = p(u,zu) truncated to u-degree at most 2j-2
        m = the largest power of u dividing \=p
        M = j*(j-1)/2;
        if (j > m+1) then M = M - (j-m)*(j-m-1)/2;
        return M
)
\enddMaccode

We implemented in Macaulay2 these algorithms
for computing the instanton widths and heights.
The computation of height of course only takes a negligible amount of time,
and the computation of widths takes a few seconds. 
For example,
 \Maccode iwidth(x^4-x*y^5,4) \endMaccode
finishes in  a Linux workstation in 17.07 seconds,
and \Maccode iwidth(x^4-x^2*y^3 -x^3*y^5-y^8,4) \endMaccode
finishes in 32.02 seconds.

%
%
%
%

\vskip3ex
\noindent
{\bf Acknowledgment.}
We thank Michael Stillman for help with Macaulay2 code.
We also thank the referees for suggesting several 
improvements to the exposition.


\bigskip
\leftline{\bf References}
\bigskip

\font\eightrm=cmr8 \def\rm{\fam0\eightrm}
\font\fiverm=cmr5 \def\rm{\fam0\eightrm}
\font\eightit=cmti8 \def\it{\fam\itfam\eightit}
\font\eightbf=cmbx8 \def\bf{\fam\bffam\eightbf}
\font\eighttt=cmtt8 \def\tt{\fam\ttfam\eighttt}
\textfont0=\eightrm \scriptfont0=\fiverm
\rm
\baselineskip=9.9pt
\parindent=3.6em

\item{[Ar]} Arnold, V. I., {\it
Singularity Theory, Selected Papers, London Math. Society Lecture 
Note Series {\bf 53} (1981)}

\item{[BG]}
Ballico, E., and Gasparim, E.,
{\it Numerical invariants for bundles on blow-ups,
Proc. Amer. Math. Soc., {\bf 130} (2002), n.1, 23--32.}

\item{[G1]} Gasparim, E.,
{\it Chern classes of bundles on blown-up surfaces, Comm. Algebra,
{\bf 28} (2000), n. 10, 4919--4926.}

\item{[G2]} Gasparim, E.,
{\it Rank two bundles on the blow up of ${\sbbC}^2$,
J. Algebra,
{\bf 199} (1998),581--590.}

\item{[G3]} Gasparim, E., {\it Two applications of instanton numbers,
Isaac Newton Institute preprint (2002),}
{\tt arXiv:math.AG/0207074.}

\item{[Ha]} Hartshorne, R.,
{\it Algebraic Geometry. 
 Graduate Texts in Mathematics} {\bf 56} Springer Verlag (1977)

\item{[Sa]}
Saito, K.,{\it  Quasihomogene isolierte Singularit\"aten von
 Hyperfl\"achen.   Invent. Math. {\bf 14}  (1971), 123--142.}

\item{[GS]}
Grayson, D., and Stillman, M., 
{\it Macaulay2,
A system for computation in algebraic geometry and commutative algebra}
(1996)
available via anonymous {\tt ftp} from {\tt math.uiuc.edu}.

\bigskip

\noindent
Department of Mathematical Sciences,
New Mexico State University,
Las Cruces, NM 88003-8001, USA, {\tt gasparim@nmsu.edu}, {\tt iswanson@nmsu.edu}.

\end